\documentclass[12pt]{article}
\usepackage{latexsym}
\usepackage{epsfig}
\usepackage{amsfonts,amssymb,amsmath}

\title{Localisable moving average symmetric\\
stable and multistable processes \\
}
\author{K.J. Falconer \\
\small{{\em Mathematical Institute,
University of St~Andrews, North Haugh, St~Andrews,}} \\
\small{{\em Fife, KY16~9SS, Scotland }} \\
R. Le Gu\'evel\\
\small{{\em Universit\'{e} Nantes, Laboratoire de Math\'{e}matiques 
Jean Leray UMR CNRS 6629}}\\
\small{{\em2 Rue de la Houssini\`{e}re - BP 92208 - F-44322 Nantes Cedex 3, France}} \\
\small{ and } \\
J. L\'{e}vy V\'{e}hel \\
\small{{\em Projet Fractales, INRIA Rocquencourt, 78153 Le Chesnay
Cedex, France}}}
\newcommand\bbbr{\mathbb{R}} % Real numbers
\newcommand\bbbz{\mathbb{Z}}
\newtheorem{theo}{Theorem}

\newtheorem{prop}[theo]{Proposition}

\newtheorem{cor}[theo]{Corollary}
\newtheorem{exa}[theo]{Example}

\newcommand\E{\mbox{\sf E}}
\newcommand\fdd{\stackrel{{\rm fdd}}{\rightarrow}}

\newcommand{\one}{\ifmmode {\sf 1}\hspace{-.26em}{\sf
l}\hspace{-.35em}{\sf \_} \else ${\sf 1}\hspace{-.26em}{\sf
l}\hspace{-.35em}{\sf \_}$ \fi}

\newcommand\varep{\varepsilon}

\newcommand\X{{\sf X}}
\newcommand\Y{{\sf Y}}
\newcommand{\1}{{\bf 1}}

\renewcommand{\Box}{\mbox{\rule{1ex}{1ex}}}

\begin{document}

\maketitle

\begin{abstract}
\noindent We study a particular class of moving average processes which possess a property called {\it localisability}. This means that, at any given point, they
admit a ``tangent process'', in a suitable sense. We give general conditions 
on the kernel $g$ defining the moving average which ensures that the process is localisable and we characterize the nature of the associated tangent processes.
Examples include the reverse Ornstein-Uhlenbeck process and the 
multistable reverse Ornstein-Uhlenbeck process. In the latter case, the 
tangent process is, at each time $t$, a L\'{e}vy stable motion with stability index
possibly varying with $t$. We also consider the problem of path synthesis, for which we
give both theoretical results and numerical simulations.
    
\end{abstract}
\section{Introduction and background}
\setcounter{equation}{0}
\setcounter{theo}{0}

In this work, we study moving average processes which are {\it localisable}. Loosely speaking, this means that they have a well-defined local form: at each point, they are ``tangent'' to a given stochastic process.

Localisable processes are useful 
both in theory and in practical applications. Indeed, they provide
an easy way to control important  
local properties such as the local H\"{o}lder 
regularity or the jump intensity. In the first case, one speaks of
{\it multifractional} processes, and in the second one, of {\it multistable}
processes. Such processes provide fine models for real world phenomena including
natural terrains, TCP traffic, financial data, EEG or highly textured images.

Formally, a process $Y(t)$ defined on $\bbbr$ (or a subinterval of $\bbbr$) is $h$-{\em localisable} at $u$ if $\lim_{r \to 0^+}r^{-h} (Y(u+rt) -Y(u))$ exists  as a non-trivial process in $t$ for some $h>0$, where the convergence is in finite dimensional distributions, see \cite{Fal5,Fal6}. When convergence occurs in distribution with respect to the appropriate metric on $C(\bbbr)$ (the space of continuous
functions on $\bbbr$) or on $D(\bbbr)$ (the space of {\it c\`{a}dl\`{a}g} functions on $\bbbr$, that is functions which are continuous on the right and
have left limits at all $t \in \bbbr$), we say that $Y$ is {\it strongly localisable}. The limit, denoted by $Y_{u}'=\{Y_{u}'(t): t\in \bbbr\}$, is called 
the {\it local form} or {\it  tangent process} of $Y$ at $u$ and will in general 
vary with $u$. A closely related notion is that of locally asymptotically
self-similar processes, which are for instance described in \cite{SC}.

The simplest localisable processes are self-similar processes with stationary 
increments (sssi processes); it is not hard to show that an sssi process $Y$ is 
localisable at all $u$ with local form $Y_{u}'=Y$. Furthermore, an sssi process $Y$ is strongly localisable if it has a version in $C(\bbbr)$ or $D(\bbbr)$.

In \cite{FLV}, processes  with prescribed local form are constructed by ``gluing together'' known localisable processes in the following way:
let $U$ be an interval with $u$ an interior point.  Let 
$\{X(t,v) : (t,v) \in U\times  U\}$ 
be a random field and  
let $Y$ be the diagonal process
$Y=\{X(t,t) : t\in U\}$. 
In order for $Y$ and $X(\cdot,u)$ to have the same local
forms at $u$, that is
$Y_{u}'(\cdot) = X_{u}'(\cdot,u)$ where $X_{u}'(\cdot,u)$ is the local
form of $X(\cdot,u)$ at $u$, we require
\begin{equation}
\frac{X(u+rt,u+rt) - X(u,u)}{r^{h}} \fdd X_{u}'(t,u) \label{locform}
\end{equation}
as $r \searrow 0$.  

This approach allows easy construction of localisable processes
from ``elementary pieces'' which are known to be themselves localisable. In
particular, it applies in a straightforward way to processes $X(t,v)$ such 
that $X(\cdot,v)$ is sssi for each $v$.

In this work we shall study a rather different way of obtaining localisable processes. 
Instead of basing our constructions on existing sssi processes we will
consider moving average processes which, as we shall see, provide a 
new class of localisable $\alpha$-stable processes.

The remainder of this paper is organized as follows: in section \ref{wl}, we
give general conditions on the kernel $g$ defining the moving average process
to ensure (strong) localisability. Section \ref{w2} specializes these conditions 
to cases where explicit forms for the tangent process may be given, and 
presents some examples. In section \ref{ms}, we deal with multistable moving
average processes, which generalize moving average stable processes by letting
the stability index vary over time. Finally, section \ref{ps} considers numerical
aspects: for applications, it is desirable to synthesize paths of these
processes. Using the approach developed in \cite{ST3}, we first explain how to build
traces of arbitrary moving average stable processes. In the case where the processes
are localisable, we then give error bounds between
the numerical and theoretical paths. Under mild additional assumptions, 
an `optimal' choice of the parameters defining the
synthesis method is derived. Finally, traces obtained
from numerical experiments are displayed. 
%222222222222

\section{Localisability of stable moving average processes}\label{wl}
\setcounter{equation}{0}
\setcounter{theo}{0}

Recall that a process $\{X(t):t\in T\}$, where $T$ is
a subinterval of $\bbbr$,
is called $\alpha$-{\it stable} $(0<\alpha \leq 2)$ 
if all its finite-dimensional distributions are
$\alpha$-stable, see the encyclopaedic work on stable 
processes \cite{ST}. $2$-stable processes are just Gaussian 
processes. 

Many stable processes admit a stochastic integral representation.
Write $S_{\alpha}(\sigma,\beta,\mu)$ for the $\alpha$-stable distribution
with scale parameter $\sigma$, skewness $\beta$ and shift-parameter 
$\mu$; we will assume throughout that $\beta = \mu=0$.  Let $(E,{\cal E},m)$ be a
sigma-finite measure space  ($m$ will be Lebesgue measure in our
examples).  
Taking $m$  as the control measure, this defines an
$\alpha$-stable random measure $M$ on $E$ such that
for $A\in {\cal E}$ we have that 
$M(A) \sim S_{\alpha}(m(A)^{1/\alpha},
0,0)$ (since $\beta =0$, the process is
symmetric).

Let
$${\cal F}_{\alpha}\equiv {\cal F}_{\alpha}(E,{\cal E}, m) 
= \{ f: f \mbox{ is measurable and } \|f\|_{\alpha} < \infty\},$$
where $\|\,\|_{\alpha}$ is the quasinorm (or norm if $1<\alpha \leq
2$) given by 
\begin{equation}
\|f\|_{\alpha}     
= \left(\int_E |f(x)|^{\alpha}m(dx)\right)^{1/\alpha}. 
    \label{normdef}
\end{equation}
The stochastic
integral of $f\in  {\cal F}_{\alpha}(E,{\cal E}, m)$ with respect to
$M$ then exists  \cite[Chapter 3]{ST} with
\begin{equation}
\int_E f(x)M(dx)\sim
S_{\alpha}(\sigma_{f},0,0),\label{alint}
\end{equation}
where
$\sigma_{f}=\|f\|_{\alpha}$.

We will be concerned with a special kind of stable processes that are 
stationary and may be
expressed as moving average stochastic integrals in the following way:
\begin{equation}
Y(t) = \int g(t-x)M(dx) \quad (t \in \bbbr),\label{ma}
\end{equation}
where $g \in {\cal F}_{\alpha}$ is sometimes called the {\it kernel} of $Y$.

Such processes are considered in several areas 
(e.g. linear time-invariant systems) 
and it is of interest to know under what conditions they are localisable. 
A sufficient condition is provided by the following proposition.

\begin{prop}\label{prop1}
Let $0<\alpha\leq 2$ and let $M$ be a symmetric $\alpha$-stable
measure on $\bbbr$ with control measure Lebesgue measure ${\cal L}$.
Let $g \in {\cal F}_{\alpha}$ and  
let $Y$ be the moving average process
\begin{equation*}
Y(t) = \int g(t-x)M(dx) \quad (t \in \bbbr).
\end{equation*}
Suppose that there exist jointly measurable functions $h(t,.) \in 
{\cal F}_{\alpha}$ such that 
\begin{equation}
\lim_{r\to 0^+} \int\left|\frac{g(r(t+z))-g(rz)}{r^{\gamma}} -
h(t,z)\right|^{\alpha}dz=0 \label{malim}
\end{equation}
 for all $t\in \bbbr$, where $\gamma+1/\alpha >0$.  Then $Y$ is
$(\gamma+1/\alpha)$-localisable with local form
$Y'_{u}= \{\int h(t,z)M(dz): t\in \bbbr\}$ at all $u\in \bbbr$.  
\end{prop}

\noindent{\it Proof.} Using stationarity followed by a change of
variable $z=-x/r$ and the self-similarity of $M$,
\begin{eqnarray*}
Y(u+rt)-Y(u) &=& Y(rt)-Y(0) \\
&=& \int (g(rt-x)-g(-x))M(dx)\\
&=& r^{1/\alpha}\int (g(r(t+z))-g(rz))M(dz),
\end{eqnarray*}
where equalities are in finite dimensional distributions.
Thus 
$$\frac{Y(u+rt)-Y(u)}{r^{\gamma+1/\alpha}} -\int h(t,z)M(dz)
= \int \left(\frac{g(r(t+z))-g(rz)}{r^{\gamma}} -
h(t,z)\right)M(dz).$$
By \cite[proposition 3.5.1]{ST} and (\ref{malim}),
$r^{-\gamma-1/\alpha}(Y(u+rt)-Y(u))\to 
\int h(t,z)M(dz)$ in probability and thus in finite dimensional
distributions.  
\Box

\medskip

A particular instance of (\ref{ma}) is the reverse Ornstein-Uhlenbeck
process, see \cite[Section 3.6]{ST}. This process provides a straightforward
application of proposition \ref{prop1}.

\begin{prop}\label{locOU} (Reverse Ornstein-Uhlenbeck
process)
Let $\lambda>0$ and $1<\alpha \leq 2$ and let  $M$ be an $\alpha$-stable
measure on $\bbbr$ with control measure ${\cal L}$.
The stationary process   
\begin{equation*}
Y(t) = \int_{t}^{\infty} \exp(-\lambda (x-t))M(dx) \quad (t \in \bbbr)
\end{equation*}
has a version in $D(\bbbr)$ that is $1/\alpha$-localisable
at all $u \in \bbbr$ with $Y_{u}' =L_{\alpha}$, where 
$L_{\alpha} (t) := \int_{0}^{t} M(dz)$
is $\alpha$-stable L\'{e}vy motion.
\end{prop}

\noindent{\it Proof.} 
The process $Y$ is a moving average process 
that may be written in the form (\ref{ma})
with $g(x) = \exp(\lambda x) \1_{(-\infty, 0]}(x)$. It is
easily verified using the dominated convergence theorem 
that $g$ satisfies (\ref{malim}) with $\gamma = 0$ and
$h(t,z) = -\1_{(-t,0]}(z)$ for $t\ge 0$ and
$h(t,z)= - \1_{(0,-t]}(z)$ for $t<0$, 
so proposition \ref{prop1} gives the conclusion
with  $Y_{u}'(t) = -M ([-t,0]) = L_{\alpha}(t)$ for $t\ge 0$
and a similar formula for $t<0$.  
%Since $Y(t) \in D(\bbbr)$ this holds in the strong sense.  
\Box

%\section{Strong Localisability}\label{sl}
%\setcounter{equation}{0}
%\setcounter{theo}{0}

Proposition \ref{prop1} gives a condition on the kernel ensuring localisability. With an additional constraint we can get strong localisability. First we need the following proposition on continuity.

\begin{prop}\label{sl1}
Let $0 < \alpha < 2$, $g \in \mathcal{F}_{\alpha}$ and let $M$ be an $\alpha$-stable symmetric random measure on $\mathbb{R}$ with control measure $\mathcal{L}$. Consider the moving average process defined by (\ref{ma}).
Suppose that $g$ satisfies, for all sufficiently small $h$,
$$\int \left|g(h-x)-g(-x)\right|^{\alpha} dx \leq c |h|^\lambda,$$
where $c>0$ and $\lambda >1$. Then $Y$ has a continuous version which satisfies a $\theta$-H\"{o}lder condition for all $\theta < (\lambda-1)/\alpha$.
\end{prop}

\medskip

\noindent{\it Proof.}
By stationarity,
\begin{eqnarray*}
	Y(t)-Y(t') & = & Y(t-t') - Y(0) \\
	& = &  \int \left(g(t-t'-x)-g(-x)\right) M(dx).
\end{eqnarray*}
So for $0<p<\alpha$
\begin{eqnarray*}
	\E|Y(t)-Y(t')|^p & \leq & c_1 \left(\int \left|g(t-t'-x)-g(-x)\right|^\alpha dx\right)^{p/\alpha} \\
	& \leq &  c_2 |t-t'|^{\lambda p/\alpha}.
\end{eqnarray*}
The result then follows from the Kolmogorov criterion by taking $p$ arbitrarily close to $\alpha$.
\Box

\begin{prop}\label{sl2}
With the same notation and assumptions as in proposition \ref{prop1}, suppose that in 
addition that $g$ satisfies, for all sufficiently small $h$,
\begin{equation}\label{assumh}
\int \left|g(h-x)-g(-x)\right|^{\alpha} dx \leq c |h|^{\alpha\gamma +1},
\end{equation}
where $c>0$ and $\gamma >0$. Then $Y$ has a version in $C(\bbbr)$ that is $(\gamma + 1/\alpha)$-strongly 
localisable with $Y'_u = \{\int h(t,z)M(dz): t\in \bbbr\}$ at all $u\in \bbbr$.
\end{prop}

\medskip

\noindent{\it Proof.}
By proposition \ref{sl1}, $Y(t)$ has a continuous version and so 
$Z_r(t) : =r^{-(\gamma + 1/\alpha)}(Y(rt)-Y(0))$ also has a continuous version. Thus, 
for $0<p<\alpha$, by stationarity and setting  $h = r|t-t'|$ sufficiently small,
\begin{eqnarray*}
\E |Z_r(t) - Z_r(t')|^p & = & \E |Z_r(t-t') - Z_r(0)|^p  \\
& = & c_1\left(\int \left| \frac{g(r(t-t')-x)- g(-x)}{r^{\gamma + 1/\alpha}}\right|^{\alpha} dx \right)^{p/\alpha} \\
& = & c_1\left(\int \frac{|g(h-x)- g(-x)|^{\alpha}}{h^{\gamma \alpha +1}} 
dx \right)^{p/\alpha} |t-t'|^{(\gamma \alpha +1)p/\alpha} \\
& \leq & c_2 |t-t'|^{(\gamma \alpha +1) p/\alpha},
\end{eqnarray*}
provided $|t-t'|$ is sufficiently small, using 
(\ref{assumh}) in the last step. We may choose $p$ sufficiently close to 
$\alpha$ so that $(\gamma \alpha +1)p/\alpha >1$.
By a Corollary to Kolmogorov's criterion (see e.g. \cite[Theorem 85.5]{RW})
the measures on $C(\bbbr)$ underlying the processes $Z_r$ are conditionally compact. 
Thus convergence in finite dimensional distributions of $Z_r$ 
to $Y'_u$ as $r \searrow 0$ implies the convergence in distribution (with $Y'_u$ necessarily having a 
continuous version). Together with localisability which follows from proposition \ref{prop1} this gives strong localisability.
\Box

Note that the reverse Ornstein-Uhlenbeck process 
is a stationary Markov process which has a version
in $D(\bbbr)$ see \cite[Remark 17.3]{Sa}. It also satisfies (\ref{assumh})
for $\alpha \geq 1$ with $\gamma=0$. However, we cannot deduce that it is strongly 
localisable since proposition \ref{sl2} is only valid for $\gamma>0$. The case
$\gamma=0$ would be interesting to deal with, but is much harder and would require different techniques.

%3333333333333333
 
\section{Sufficient conditions for localisability and examples}\label{w2}
\setcounter{equation}{0}
\setcounter{theo}{0}

For the reverse Ornstein-Uhlenbeck process, it was straightforward to check
the conditions of proposition \ref{prop1}. In general, however, it is not easy to guess which kind of functions $g$ in 
$\mathcal{F}_{\alpha}$ will satisfy
(\ref{malim}).
In this section we will find simple practical conditions 
ensuring this.

Recall that the following process is called \textit{linear fractional $\alpha$-stable motion}:
\begin{displaymath}
L_{\alpha,H,b^{+},b^{-}} (t) = \int_{-\infty}^{\infty} f_{\alpha,H}(b^{+},b^{-},t,x) M(dx),
\end{displaymath}
where $t \in \mathbb{R}$, $b^{+},b^{-} \in \mathbb{R}$, and
\begin{align}
f_{\alpha,H}(b^{+},b^{-},t,x) = b^{+} & \left( (t-x)_{+}^{H - 1/ \alpha} - (-x)_{+} ^{H - 1/ \alpha} \right) \nonumber
\\
& + b^{-} \left( (t-x)_{-} ^{H- 1/ \alpha} - (-x)_{-} ^{H- 1/ \alpha} \right), \label{lfsm}
\end{align}
where $M$ is a symmetric $\alpha$-stable random measure $( 0 < \alpha < 2)$ with control measure Lebesgue measure. Being sssi, $L_{\alpha,H,b^{+},b^{-}}$ is
localisable. In addition, it is strongly localisable when $H > 1/\alpha$, since
its paths then belong to $C(\bbbr)$.

Recall also that the process
\begin{equation}
 L_{\alpha} (t) = \int_{0}^{t} M(dx)\label{slm}
\end{equation}
is  $\alpha$-{\em stable L\'evy motion} and the process
\begin{equation}
 Z_{\alpha} (t) = \int_{-\infty}^{+\infty} (\ln |t-x| - \ln |x| ) M(dx)\label{logfsm}
\end{equation}
is called {\em log-fractional stable motion}.

We are now ready to describe easy-to-check conditions that ensure that propositions \ref{prop1} and \ref{sl2} apply.

\begin{prop}\label{prop1.2}
Let $0 < \alpha \leq 2$, $g \in \mathcal{F}_{\alpha}$ and $M$ be an $\alpha$-stable symmetric random measure on $\mathbb{R}$ with control measure $\mathcal{L}$. Let $Y$ be 
the moving average process 
\begin{displaymath}
Y(t) = \int g(t-x) M(dx) \hspace{0.5cm} (t \in \mathbb{R}).
\end{displaymath}
If there exist $c_{0}^{+}, c_{0}^{-}, \gamma, a, c, \eta \in \mathbb{R}$ with $c >0$, $\eta > 0$ and $ 0< \gamma+ 1/\alpha < a \leq 1$  such that 
\begin{displaymath}
\frac{g(r)}{r^{\gamma}} \rightarrow c_{0}^{+} \textrm{ and } \frac{g(-r)}{r^{\gamma}} \rightarrow c_{0}^{-}
\end{displaymath}
as $r \searrow 0$ and
\begin{equation}\label{2ml}
\left| g(u+h) - g(u) \right| \leq c |h|^{a}|u|^{\gamma - a} \quad
(u \in \mathbb{R}, |h| <\eta), 
\end{equation}
then $Y$ is $(\gamma + 1/ \alpha)$-localisable at all $u \in \mathbb{R}$ with local form
\newline
\begin{displaymath}
\begin{array}{ll}
(a) & Y'_u = L_{\alpha, \gamma+1/\alpha,c_0^+,c_0^-}\quad \textrm{if } \gamma \neq 0,\\
 & \\
(b) & Y'_u = ( c_{0}^{+}-c_{0}^{-}) L_{\alpha} \quad \textrm{if } \gamma = 0.\\
\end{array}
\end{displaymath}
If, in addition, $\gamma>0$ and $0<\alpha<2$ then $Y$ has a version in $C(\bbbr)$ and is strongly localisable.
\end{prop}

Note that condition (\ref{2ml}) on the increments of $g$ may be interpreted as a 2-microlocal condition, namely that $g$ belongs to the global 2-microlocal space $C_0^{\gamma,a - \gamma}$, see \cite{LVS}. Remark also that, in order for this condition to be satisfied by non-trivial functions $g$, one needs $a\leq 1$, which in turns implies that $\gamma \leq 1 - 1/\alpha$ and $a - \gamma \in (1/\alpha, 1-\gamma]$.

\bigskip

\noindent{\it Proof} 
(a) We have
\begin{align*}
\frac{g(r(t+z))-g(rz)}{r^{\gamma}} = \frac{g(r|t+z|)}{(r|t+z|)^{\gamma}}&|t+z|^{\gamma}\mathbf{1}_{\{ t+z \geq 0\}} + \frac{g(-r|t+z|)}{(r|t+z|)^{\gamma}}|t+z|^{\gamma}\mathbf{1}_{\{ t+z < 0\}}
\\& - \frac{g(r|z|)}{(r|z|)^{\gamma}}|z|^{\gamma}\mathbf{1}_{\{ z \geq 0\}} - \frac{g(-r|z|)}{(r|z|)^{\gamma}}|z|^{\gamma}\mathbf{1}_{\{ z < 0\}}.
\end{align*}
As $r \to 0$,
\begin{eqnarray}
\frac{g(r(t+z))-g(rz)}{r^{\gamma}} & \rightarrow & c_{0}^{+}|t+z|^{\gamma}\mathbf{1}_{\{ t+z \geq 0\}} + c_{0}^{-}|t+z|^{\gamma}\mathbf{1}_{\{ t+z < 0\}} - c_{0}^{+}|z|^{\gamma}\mathbf{1}_{\{ z \geq 0\}} - c_{0}^{-}|z|^{\gamma}\mathbf{1}_{\{ z < 0\}}\nonumber\\
& = & c_{0}^{+} (t+z)_{+}^{\gamma} - c_{0}^{+} (z)_{+}^{\gamma} +c_{0}^{-}(t+z)_{-}^{\gamma} - c_{0}^{-}(z)_{-}^{\gamma} \label{limf}\\
&=&f_{ \alpha, \gamma + 1/ \alpha} (c_{0}^{+},c_{0}^{-},t,-z).\nonumber
\end{eqnarray}
To get convergence in $L^{\alpha}$ we use the  dominated convergence theorem.
Fix $\epsilon > 0$ and $m > 0$ such that for all $0 < u < \epsilon$,
\begin{displaymath}
 \left| \frac{g(u)}{u^{\gamma}} - c_{0}^{+} \right| \leq m \textrm{  and  }  \left| \frac{g(-u)}{u^{\gamma}} - c_{0}^{-} \right| \leq m.
\end{displaymath}
For fixed $t$ write $f_r(z) = r^{-\gamma}(g(r(t+z))-g(rz))$. If $t=0$, $$f_r(z)=f_{ \alpha, \gamma + 1/ \alpha} (c_{0}^{+},c_{0}^{-},t,-z)= 0,$$ thus $f_r \to f_{ \alpha, \gamma + 1/ \alpha} (c_{0}^{+},c_{0}^{-},t,-.)$ which belongs to $L^{\alpha}$. Assume now $t \in \bbbr^*$. There is a constant $m_1$ such that  
$|f_r(z)|^{\alpha} \leq  m_1(1+|z|^{\gamma}+|t+z|^{\gamma})^{\alpha}$ for all $|r| \leq \frac{\epsilon}{1+|t|}$ and $|z| \leq 1$.
From (\ref{2ml})
$$|f_r(z)|^{\alpha} \leq \left( \frac{|rt|^{a}|rz|^{\gamma -a}}{ |r|^{\gamma}}\right)^{\alpha}
\leq|t|^{a\alpha}|z|^{(\gamma -a)\alpha}$$
 for $|r| < \eta/|t|$, so, as $( \gamma-a )\alpha <-1$ and $ \gamma \alpha > -1$,
 $$ \int_{|z| \leq1}m_1(1+|z|^{\gamma}+|t+z|^{\gamma})^{\alpha} dz +     \int_{|z| >1}   |z|^{(\gamma -a)\alpha} dz < \infty.$$
 
Since also $f_{ \alpha, \gamma + 1/ \alpha} (c_{0}^{+},c_{0}^{-},t,-z) \in L^{\alpha}$, the dominated convergence theorem implies that 
$f_r(z) \to f_{ \alpha, \gamma + 1/ \alpha} (c_{0}^{+},c_{0}^{-},t,-z)$ in $L^{\alpha}$. 
The conclusion in case (a) follows from proposition \ref{prop1}, (\ref{lfsm}), and noting that $M$ is a symmetric $\alpha$-stable measure.

\medskip
(b)
In this case the limit (\ref{limf}) is
$$
\frac{g(r(t+z))-g(rz)}{r^{\gamma}}
\to \left\{
\begin{array}{rcl}
 (c_{0}^{+} - c_{0}^{-}) \mathbf{1}_{[0,t]}(-z) &\mbox{ if }  & t\geq 0   \\
 -(c_{0}^{+} - c_{0}^{-}) \mathbf{1}_{[t,0]}(-z) &  \mbox{ if }&  t<0   \\
\end{array} \right..
$$
Dominated convergence follows in the same way as in case (a) so the conclusion follows from proposition \ref{prop1} and (\ref{slm}).

Moving to  strong localisability,
for $h$ small enough,
$$\int_{|x| \leq 3|h|}|g(h-x)-g(-x)|^{\alpha}dx \leq c_1 \int_{|x| \leq 3|h|} |h|^{\alpha\gamma} dx \leq c_2 |h|^{\alpha\gamma+1},$$
and 
\begin{eqnarray*}
\int_{|x| \geq 3|h|}|g(h-x)-g(-x)|^{\alpha}dx & \leq & c_1 |h|^{a\alpha}\int_{3|h|}^{\infty} |x|^{(\gamma-a)\alpha} dx \\
& \leq & c_2 |h|^{a\alpha}|h|^{1+(\gamma-a)\alpha} \\
& = & c_2 |h|^{\alpha\gamma+1}
\end{eqnarray*}
and the conclusion follows from propositions \ref{sl1} and \ref{sl2}. \Box
\medskip

We now give an alternative condition for localisability in terms of Fourier transforms. Note that 
the Fourier transform $\widehat{f}_{\alpha,H} (b^{+},b^{-},t,\xi)$ of 
$f_{\alpha,H}(b^{+},b^{-},t,.)$ is given by 
\begin{align*}
\widehat{f}_{\alpha,H} (b^{+},b^{-},t,\xi)=  & \Gamma(H +1 -1/\alpha) 
\frac{e^{-i\xi t}-1}{|\xi|^{H +1 -1/\alpha}} \\
&\times \left[  b^{+}\exp\Big(\frac{i \pi}{2} \textrm{sgn}(\xi) (H +1 -1/\alpha)\Big)  
+  b^{-}\exp\Big(-\frac{i \pi}{2} \textrm{sgn}(\xi)( H +1 -1/\alpha\Big)\right].
\end{align*}

\begin{prop}\label{prop1.2A}
Let $1 \leq \alpha \leq 2$, and $Y$ be defined by (\ref{ma}).
If there exist $l = l_1 + il_2 \in \mathbb{C}^{*}$, $\gamma \in (-\frac{1}{\alpha},1-\frac{1}{\alpha})$, $a \in (0,1-(\gamma + \frac{1}{ \alpha}))$ and $K \in L^p(\mathbb{R})$ with $ p \in [1 , 1/(\gamma + \frac{1}{\alpha} + a)) $,  such that for almost all $\xi > 0$,
\begin{equation}\label{condufourier}
 \xi^{\gamma +1} \widehat{g}(\xi) = l + \frac{1}{\xi^a}\widehat{K}(\xi),
\end{equation}
then $Y$ is $(\gamma + 1/ \alpha)$-localisable at all $u \in \mathbb{R}$ with local form
\newline
\begin{displaymath}
\begin{array}{ll}
(a) & Y'_u = L_{\alpha, \gamma+1/\alpha,b^+,b^-}\quad  \textrm{if } \gamma \neq 0\\
 & \\
(b) & Y'_u = \frac{1}{\pi} l_1 Z_{\alpha} +l_2 L_{\alpha} \quad \textrm{if } \gamma = 0,\\
\end{array}
\end{displaymath}
where 
\begin{displaymath}
b^+ = \frac{1}{2 \Gamma(\gamma +1)} \left( \frac{l_1}{\cos (\pi (\gamma +1)/2)} - \frac{l_2}{\sin  (\pi (\gamma +1)/2)} \right),
\end{displaymath}
\begin{displaymath}
b^- = \frac{1}{2 \Gamma(\gamma +1)} \left( \frac{l_1}{\cos (\pi (\gamma +1)/2)} + \frac{l_2}{\sin  (\pi (\gamma +1)/2)} \right).
\end{displaymath}
\end{prop}

\noindent{\em Proof} (a)
First note that, with $b^+$ and $b^-$ as above, we have, for $z \neq 0$,
\begin{displaymath}
\widehat{f}_{\alpha,\gamma + 1/ \alpha} (b^+,b^-,t,\xi) = \frac{e^{- i \xi t}-1}{|\xi|^{\gamma+1}} (\bar{l} \mathbf{1}_{\xi > 0} + l\mathbf{1}_{\xi < 0} ).
\end{displaymath}
Set $ f_r(z) = r^{-\gamma}(g(r(t+z))-g(rz))$. Then $f_r \in \mathcal{F}_{\alpha}$ and
\begin{displaymath}
\widehat{f}_r(\xi) = \frac{e^{i \xi t}-1}{r^{\gamma +1}} \widehat{g}\big(\frac{\xi}{r}\big).
\end{displaymath}
With $\alpha'$ such that $\frac{1}{\alpha}+\frac{1}{\alpha'}=1$ we have $\widehat{f}_r \in \mathcal{F}_{\alpha'}$ and $\widehat{f}_{\alpha,\gamma + 1/ \alpha} (b^+,b^-,t,\xi) \in \mathcal{F}_{\alpha'}$. We now show that $\Vert f_r - f_{\alpha,\gamma + 1/ \alpha} (b^+,b^-,t,-\cdot)\Vert_{\alpha} \rightarrow 0$ when $r \rightarrow 0$.  Note that (\ref{condufourier}) implies that for $\xi<0$
\begin{equation*}
 |\xi|^{\gamma +1} \widehat{g}(\xi) = \bar{l} + \frac{1}{|\xi|^a}\widehat K(\xi).
\end{equation*}
Writing $\widehat{f}(\xi) = \widehat{f}_{\alpha,\gamma + 1/ \alpha} (b^+,b^-,t,-\xi)$,
for almost all $\xi \in \mathbb{R}$ 
\begin{eqnarray*}
\widehat{f}_r(\xi) - \widehat{f}(\xi)  & = &  \frac{(e^{i \xi t}-1)}{|\xi|^{\gamma+1}} \left( \Big( \frac{|\xi|^{\gamma+1}}{r^{\gamma+1}} \widehat{g}\big(\frac{\xi}{r}\big) - l\Big) \mathbf{1}_{\xi > 0} + \Big(\frac{|\xi|^{\gamma+1}}{r^{\gamma+1}} \widehat{g}\big(\frac{\xi}{r}\big) - \bar{l} \Big) \mathbf{1}_{\xi < 0} \right)\\
& = & \frac{(e^{i \xi t}-1)}{|\xi|^{\gamma+1}} \frac{r^a}{|\xi|^a} K\big(\frac{\xi}{r}\big) \\
& = & r^a \frac{(e^{i \xi t}-1)}{|\xi|^{\gamma+1+a}}K\big(\frac{\xi}{r}\big).\\
\end{eqnarray*}
Let $H_r(u) = K(ru)$. Then $\widehat{H}_r(\xi)= \frac{1}{r}\widehat{K}(\frac{\xi}{r})$ and we may write for $a + \gamma \neq 0$
\begin{eqnarray}\label{goodeq}
 \widehat{f}_r(\xi) - \widehat{f}(\xi) & = & r^{a+1} \widehat{f}_{\alpha,\gamma + 1/ \alpha + a} (b,b,t,-\xi)\widehat{H_r}(\xi),
\end{eqnarray}
where $b  = 1/(2 \Gamma(\gamma+a+1) \cos(\pi(\gamma+a+1)/2))$. 

It is easy to verify that $f_{\alpha,\gamma + a+ 1/ \alpha} (b^+,b^-,t,-\cdot) \in L_{\beta}$ for all $\beta > 1/(1-\gamma - a)$. By the conditions on $\alpha$ and $p$, there exists such a $\beta$ which also satisfies $\frac{1}{\alpha}+1=\frac{1}{\beta} + \frac{1}{p}$ and in particular, $\frac1 p + \frac 1 \beta >1$. Consequently we may take the inverse Fourier transform of (\ref{goodeq}) see, for example,  \cite[Theorem 78]{Tit} to get:  
$$
 f_r(z) - f(z) = r^{a+1} f_{\alpha,\gamma + 1/ \alpha + a} (b^+,b^-,t,-.)*H_r (z)
$$
where $*$ denotes convolution.
As $\frac{1}{\alpha}+1=\frac{1}{\beta} + \frac{1}{p}$, the Hausdorff-Young inequality yields
\begin{eqnarray*}
 \Vert f_r - f_{\alpha,\gamma + 1/ \alpha} (b^+,b^-,t,-\cdot)\Vert_{\alpha} & \leq & r^{a+1} \Vert f_{\alpha,\gamma + 1/ \alpha + a} \Vert_{\beta} \Vert H_r \Vert_{p} \\
 & \leq & r^{a +1 - \frac{1}{p}} \Vert f_{\alpha,\gamma + 1/ \alpha + a} \Vert_{\beta}\Vert K \Vert_{p}. 
\end{eqnarray*}
We conclude that $f_r \rightarrow f_{\alpha,\gamma + 1/ \alpha} (b^{+},b^{-},t,-\cdot)$ in  $L^{\alpha}$. The result follows from proposition \ref{prop1}.
The case $a + \gamma = 0$ is dealt with in a similar way.

\medskip
(b)
Let $z_t$ and $l_t$ be defined by
\begin{displaymath}
l_t(x) = \left\{ 
\begin{array}{rl}
\mathbf{1}_{]0,t[}(x) & \mbox{ if }  t \geq 0\\
 -\mathbf{1}_{]t,0[}(x) & \mbox{ if }  t < 0\\
\end{array}
\right.
\end{displaymath}  
and
\begin{displaymath}
 z_t(x) = \ln|t-x| - \ln|x|.
\end{displaymath}
A straightforward computation shows that
\begin{displaymath}
 z_t(x)  =  \mbox{sgn}(-t) \lim_{\varep \rightarrow 0} \int_{|s| \geq \varep} \frac{1}{s} \mathbf{1}_{[\min(x-t,x),\max(x-t,x)]}(s) ds, 
\end{displaymath}
so that, in the space of distributions we get
\begin{displaymath} 
z_t  =  - \mbox{PV}(1/\cdot) * l_t
\end{displaymath}
where PV denotes the Cauchy principal value. Thus 
\begin{eqnarray*}
 \widehat{z_t}(\xi) &=& -\widehat{\mbox{PV}(1/\cdot)}(\xi)\widehat{l_t}(\xi)\\
 & =& -(-i \pi \mbox{sgn}(\xi))\big(-\frac{1}{i \xi} (e^{-i \xi t}-1)\big)\\
 & =& - \pi \frac{e^{-i \xi t} -1}{|\xi|}.\\
\end{eqnarray*}
With $f(z) = -\frac{1}{\pi}l_1 z_t(-z) - l_2 l_t(-z)$, we obtain 
\begin{displaymath}
 \widehat{f}(\xi) = \frac{e^{i\xi t}-1}{|\xi|}(l \mathbf{1}_{\xi > 0} + \bar{l}\mathbf{1}_{\xi < 0} ).
\end{displaymath}
As in (a) we conclude that $f_r \rightarrow f$ in  $L^{\alpha}$. proposition \ref{prop1} implies that $Y$ is $(\gamma + 1/ \alpha)$-localisable at all $u \in \mathbb{R}$ with local form 
$Y'_u = \frac{1}{\pi}l_1 Z_{\alpha} +l_2 L_{\alpha}$, since $M$ is symmetric.
\Box
\medskip

We give examples to illustrate  propositions \ref{prop1.2} and \ref{prop1.2A}.
\begin{exa}\label{extime}
 Let $ \frac{6}{5} < \alpha \leq 2$ and let $M$ be an $\alpha$-stable symmetric random measure on $\mathbb{R}$ with control measure $\mathcal{L}$. Let
 \begin{displaymath}
  g(x) =\left\{ 
\begin{array}{c c}
0 &   (x \leq 0)\\
x^{1/6} &  (0 < x \leq 1)\\
x^{-5/6} &  (x \geq 1)\\
\end{array}
\right. .
 \end{displaymath}
The stationary process defined by 
\begin{displaymath}
 Y(t) = \int g(t-x) M(dx)  \hspace{0.4cm} (t \in \mathbb{R})
\end{displaymath}
is $(1/6 + 1/ \alpha)$-strongly localisable at all $u \in \mathbb{R}$ with local form $Y'_u = L_{\alpha,1/6 + 1/\alpha,1,0}$.
\end{exa}
\medskip

\noindent{\em Proof.}
We apply proposition \ref{prop1.2}  case (a) with $\alpha \in (\frac{5}{6}, 2]$. The function $g$ satisfies the assumptions with $\gamma = \frac{1}{6}$, $c_0^+=1$, $c_0^-=0$ and $a=1$.
\Box

\medskip

To verify condition (\ref{condufourier})  of proposition  \ref{prop1.2A}, one needs to check that $g \in L^{\alpha}(\bbbr)$ and also that $\xi^{a+\gamma+1}\widehat g(\xi) - l \xi^a$ is the Fourier transform of a function in $L^p(\bbbr)$ for some $a,\gamma,p$ in the admissible ranges. For this purpose, one may for instance apply classical theorems such as in \cite[Theorems 82-84]{Tit}. We give below an example that uses a direct approach. 

\begin{exa}\label{exfrequency}
For $1\leq \alpha<2 $ let $M$ be an $\alpha$-stable symmetric random measure on $\mathbb{R}$ with control measure $\mathcal{L}$.  
Let $g$ be defined by its Fourier transform
 \begin{displaymath}
  \widehat g(\xi) =\left\{ 
\begin{array}{c c}
0 & ( |\xi| \leq 1)\\
|\xi|^{-\gamma -1} & ( |\xi| > 1)
\end{array}
\right.
 \end{displaymath}
 where  $\gamma \in (-\frac{1}{\alpha},\frac 1 2-\frac{1}{\alpha})\subseteq (-1,0)$.
Then  $g \in L^{\alpha}(\bbbr)$ and 
the moving average process
\begin{displaymath}
 Y(t) = \int g(t-x) M(dx)  \hspace{0.4cm} (t \in \mathbb{R})
\end{displaymath}
is well-defined and 
$\alpha$-localisable at all $u \in \mathbb{R}$, with local form $Y'_u = L_{\alpha,\gamma + 1/\alpha,b,b}$,
where $b=-1/(2\Gamma(\gamma+1)\cos(\pi(\gamma+1/2)))$.
\end{exa}

\noindent{\em Proof.}
Taking  $\widehat K(\xi) =|\xi|^{1/2} \mathbf{1}_{[-1,1]}(\xi)$ with $l=-1$ and $a=\frac 12$ in (\ref{condufourier}) gives $g$.
To check that $K \in L^p(\bbbr)$ for all $p>1$, note 
that $K$ is continuous 
(in fact $C^{\infty}$) and that
$|K(x)| \leq C|x|^{-1}$ for all $x$. Then $ Y(t)$
will be well-defined if $g$ is in $L^{\alpha}(\bbbr)$. To verify this,
one computes the inverse Fourier transform of $\widehat g$, to get
$ g(x)= 2(\gamma+1) |x|^{\gamma} \int_{|x|}^{\infty} |v|^{-\gamma-2}\sin v dv - 2x^{-1}\sin x$.
By proposition \ref{prop1.2A}(a), $Y$ is
$\alpha$-localisable at all $u \in \mathbb{R}$ with the local form as stated. \Box
\medskip

The approach of this example may be be used for general classes of functions $g$.

\section{Multistable moving average processes}\label{ms}
\setcounter{equation}{0}
\setcounter{theo}{0}

In \cite{FLV}, localisability is used to define {\it multistable processes}, that is processes which  at each point $t \in \bbbr$ have an 
$\alpha(t)$-stable random process as their local form, where $\alpha(t)$ is a sufficiently 
smooth function ranging in
$(0,2)$. Thus such processes  ``look locally like'' a stable 
process at each $t$ but with
differing stability indices as time evolves.

\medskip

Before we recall how this was done in \cite{FLV}, we note briefly 
that ``stable-like'' processes have been defined and studied in \cite{Neg}.
These stable-like processes are Markov jump processes,
and are, in a sense, ``localisable'', but with localisability defined 
by the requirement that they are solutions of an order $\alpha(x)$ fractional stochastic differential equations. See Theorem 2.1 in \cite{Neg}, which
shows that the local form of sample paths is considered rather than of the limiting process. Another essential difference is that stable-like processes are Markov, whereas, in general, multistable ones, as defined below, are not. In fact, formula (\ref{mmap}), where a Poisson process element Y is independent
of  t but is raised to a  power that involves t means that our processes are ``far'' from Markov.

\medskip
 
We now come back to our multistable processes. One route to defining
such processes is to rewrite 
stable integrals as countable sums over Poisson processes. We recall
briefly how this can be done, see \cite{FLV} for fuller details. Let $(E,{\cal E},m)$ be a 
$\sigma$-finite measure space and let
$\Pi$ be a Poisson process on $E \times \bbbr$ with mean measure
$m \times {\cal L}$.  
Thus $\Pi$ is a random countable subset of $E \times \bbbr$ such that, 
writing $N(A)$ for the number of points in a measurable 
$A \subset E \times \bbbr$, the random variable
$N(A)$ has a Poisson distribution of mean $(m \times {\cal L})(A)$
with $N(A_{1}),\ldots,N(A_{n})$ independent
for disjoint $A_{1},\ldots,A_{n}\subset E \times \bbbr$, see \cite{Ki}. 
In the case of constant $\alpha$, with $M$  a symmetric $\alpha$-stable
random measure on $E$ with control measure $m$,
one has, for $f \in {\cal F}_{\alpha}$  (\cite[Section 3.12]{ST}), 
\begin{equation}
\int f(x) M(dx)\, = \, c(\alpha)
\sum_{(\X,\Y) \in \Pi} f(\X) \Y^{<-1/\alpha>}\quad (0<\alpha<2),\label{ppal}
\end{equation}
where 
\begin{equation}
    c(\alpha) =
\left(2\alpha^{-1}\Gamma(1-\alpha)
\cos(\textstyle{\frac{1}{2}}\pi
\alpha)\right)^{-1/\alpha},\label{calpha}
\end{equation}
and $a^{<b>} = \mbox{sign}(a)|a|^{b}$.

Now define the random field
\begin{equation}
X(t,v)= \sum_{(\X,\Y) \in \Pi} f(t,v,\X) \Y^{<-1/\alpha(v)>}\label{msfielda}.
\end{equation}
Under certain conditions the ``diagonal'' process  $X(t,t)$ gives rise to a multistable process with varying $\alpha$ of the form
\begin{equation}
Y(t) \equiv X(t,t)= \sum_{(\X,\Y) \in \Pi} f(t,t,\X)
\Y^{<-1/\alpha(t)>}.\label{msproc}
\end{equation}
Theorem 5.2 of \cite{FLV} gives conditions on $f$ that ensure that $Y$ is localisable (or strongly localisable) with $Y'_u=X'_u(\cdot,u)$ at a given $u$, provided $X(\cdot,u)$ is itself localisable (resp. strongly localisable) at $u$. 
These conditions simplify very considerably in the moving average case, taking  $E=\bbbr$  and $m = {\cal L}$ with
$f(t,v,x) = g(x-t)$. Our next theorem restates  \cite[Theorem 5.2]{FLV} in this specific situation.

We need first to define a quasinorm on certain spaces of measurable functions on $E$. 
For $ 0<a \leq b <2$ let
$${\cal F}_{a,b}\equiv {\cal F}_{a,b}(E,{\cal E},m)
= \{f: f \mbox{ is $m$-measurable with } \|f\|_{a,b}<\infty\}$$
where
\begin{equation}
\|f\|_{a,b} = \left(\int_{E}|f(x)|^{a}m(dx)\right)^{1/a}
+  \left(\int_{E}|f(x)|^{b}m(dx)\right)^{1/b}.\label{norm}
\end{equation}

\begin{theo}\label{thmflv}{(Multistable moving average
processes)}
Let $U$ be a closed interval with $u$ an interior point. Let $\alpha : U \to (a,b) \subset (0, 2)$ satisfy 
\begin{equation*}
    |\alpha(v) - \alpha(u) | \leq k_{1} |v-u|^{\eta} \quad (v \in U)
    \end{equation*}  
where $0 < \eta \leq 1$. Let $g \in {\cal F}_{a,b}$, and define
\begin{equation}\label{mmap}
Y(t)= \sum_{(\X,\Y) \in \Pi} 
g(\X-t) \Y^{<-1/\alpha(t)>}  \quad (t \in \bbbr).
\end{equation}
Assume that $g$ satisfies
\begin{equation}
\lim_{r\to 0} \int\left|\frac{g(r(t+z))-g(rz)}{r^{\gamma}} -
h(t,z)\right|^{\alpha(u)}dz=0 \label{htz} 
\end{equation}
for jointly
measurable functions $h(t,\cdot) \in {\cal F}_{\alpha(u)}$,
where $0<\gamma +1/\alpha(u) < \eta\leq 1$.
Then $Y$ is
$(\gamma+1/\alpha(u))$-localisable at $u$ with local form
$Y'_{u}= \{\int h(t,z)M_{\alpha(u)}(dz): t\in \bbbr\}$, where $M_{\alpha(u)}$
is the symmetric $\alpha(u)$-stable measure with control measure ${\cal L}$ and skewness $0$. 

Suppose further that $\gamma>0$ and for $h$ sufficiently small 
$$\|g(h-x)-g(-x)\|_\alpha \leq c|h|^{\gamma + 1/\alpha(u)} .$$
Then $Y$ has a continuous version and is strongly
$(\gamma+1/\alpha(u))$-localisable at $u$ with local form
$Y'_{u}= \{\int h(t,z)M_{\alpha(u)}(dz): t\in \bbbr\}$
under either of the following additional conditions:

$(i)$ $0<\alpha(u)<1$ and $g$ is bounded 
%$($so if $g \in {\cal F}_{a}$ then necessarily $g \in {\cal F}_{a,b}$$)$.

$(ii)$ $1<\alpha(u)<2$ and $\alpha$ is continuously differentiable on $U$ with 
$$ |\alpha'(v) - \alpha'(w) | \leq k_{1} |v-w|^{\eta} \quad (v,w \in U),$$
where $1/\alpha(u)< \eta \leq 1$. 
\end{theo}

\noindent{\it Proof.} 
Taking
\begin{equation}
X(t,v) = \sum_{(\X,\Y) \in \Pi} g(\X-t)
\Y^{<-1/\alpha(v)>}  \quad (t,v \in \bbbr).
\end{equation} 
this theorem is essentially a restatement of \cite[Theorem 5.2]{FLV}
in the special case of $E=\bbbr$  and $m = {\cal L}$ with
$f(t,v,x) = g(x-t)$ in (\ref{msfielda}). Since  $f(t,v,x)$ no longer depends on $v$ most of the conditions  in \cite[Theorem 5.2]{FLV} are trivially satisfied and we conclude that $Y'_u=X'_u(\cdot,u)$, noting that  $X(\cdot,u)$ is $(\gamma+1/\alpha(u))$-localisable (or strongly localisable) with the local form given   
by propositions \ref{prop1} or \ref{sl2}.
\Box

\medskip

It is curious that neither cases (i) or (ii) address localisability if $\alpha(u)=1$. This goes back to the proof of [3, Theorem 5.2] where different approaches are used in the two cases. For  $\alpha(u) < 1$ the proof uses that the sum (4.8) is absolutely convergent almost surely, whereas for $\alpha(u) > 1$ we need to find a number $p$ such that $0<p<\alpha(t)$ for $t$ near $u$  with $\eta p >1$ to enable us to apply Kolmogorov's criterion to certain increments.

\medskip
 
\begin{cor}\label{mapvp}
Let $U, \alpha$ and $g$ be as in Theorem \ref{thmflv}. Then the same conclusion holds if $Y(t)$ in (\ref{mmap}) is replaced by 
$Y(t) = a(t) \sum_{(X,Y)\in \Pi} g(\X-t)\Y^{<-1/\alpha(t)>} \quad (t \in \bbbr)$,
where $a$ is a non-zero function of H\"{o}lder exponent $\eta > h$.
\end{cor}

\noindent{\it Proof.} 
This follows easily in just the same way as proposition 2.2 of \cite{FLV}.
\Box

\medskip

We may apply this theorem to get a multistable version of the reverse Ornstein-Uhlenbeck
process considered in Section \ref{wl}:

\begin{prop} (Multistable reverse Ornstein-Uhlenbeck
process)
Let $\lambda>0$ and $\alpha : \bbbr \to (1, 2)$ be continuously
differentiable.  Let
\begin{equation*}
Y(t)= \sum_{(\X,\Y) \in \Pi, \X \geq t} 
\exp(-\lambda(\X-t)) \Y^{<-1/\alpha(t)>}  \quad (t \in \bbbr).
\end{equation*}
Then $Y$ is $1/\alpha(u)$-localisable at all $u \in \bbbr$ 
with $Y'_u = c(\alpha(u))^{-1}L_{\alpha(u)}$, where $L_{\alpha}$ is 
$\alpha$-stable L\'{e}vy motion.  
\end{prop}

\noindent{\it Proof.} 
Taking $g(x) = \1_{[0,\infty)}(x)\exp(-\lambda x)$ and
$h(t,z) =  -\1_{[-t,0]}(z)$ for $t\ge 0$ (and a similar formula for $t <0$) with $\gamma = 0$, localisability follows from Theorem \ref{thmflv} with  the limit 
(\ref{htz}) being checked just as in proposition \ref{locOU}.\Box

Theorem \ref{thmflv} applies in particular to functions $g$ satisfying the conditions
of proposition \ref{prop1.2}. Thus, for instance, the moving averages of Examples
\ref{extime} and \ref{exfrequency} admit multistable versions. The
process of Example \ref{extime} is strongly $\gamma + 1/\alpha(u)$ localisable at $u$ whenever $\alpha$ verifies condition {\it (ii)}. 

%66666666666666666
\section{Path synthesis and numerical experiments}\label{ps}
\setcounter{equation}{0}
\setcounter{theo}{0}

We address here the issue of path simulation.
In the previous sections, we have considered two kinds of stochastic processes:
moving average stable ones, that are stationary, and
their multistable versions, which typically are not, nor have stationary increments.
Our simulation method for the moving average stable processes is based on that
presented in \cite{ST3}. There, the authors propose an efficient algorithm
for synthesizing paths of linear fractional stable motion. In fact, 
this algorithm really builds traces of the increments of linear fractional stable motion.
These increments form a stationary process, an essential feature for
the algorithm to work. It is straightforward to modify it to synthesize any stationary stable process which possesses an 
integral representation. In addition, we are able to obtain bounds on
the approximation error measured in the $\alpha$-norm, and thus
on the $r^{th}$ moments for $r < \alpha$, as shown below.

For non (increment) stationary processes, like multistable processes, 
a possibility would be to use the general method proposed in \cite{CLL}.
It allows to synthesize (fractional) 
fields defined by integration of a deterministic kernel with 
respect to a random infinitely divisible measure.
When the control measure is finite, the idea is to approximate the
integral with a generalized shot noise series. In this situation, a bound 
on the $L^r$ norm of the error is obtained for appropriate $r$. In the 
case of infinite control measure, one needs to deal with the points ``far from
the origin'' through a normal approximation. This second approximation
maybe controlled through Berry-Esseen bounds which lead to a convergence
in law. Thus the overall error when the control measure is infinite may only
be assessed in law, and not in the stronger $L^r$ norm. 

Although the method of \cite{CLL} may be used for the synthesis of multistable processes, we will rather take advantage here of the particular structure of our processes: being localisable, they are by definition tangent,
at each point, to an increment stationary process. Thus we may simulate
them by ``gluing'' together in a appropriate way paths of their tangent processes, which are themselves synthesized through the simpler procedure of \cite{ST3}. This approach allows in addition to control the error in the $L^r$ norm, rather than only in law (since we are in an infinite control measure case).

We briefly present in the next subsection the main ingredients of the method. 
We then give bounds estimating the errors entailed by 
the numeric approximation, in the case where the process is localisable. 
Finally, we display  graphs of localisable moving average processes obtained 
with this synthesis scheme.

\subsection{Simulation of stable moving averages}

Let $Y=\{ Y(t), t \in \bbbr \}$ be the process defined by (\ref{ma}).
To synthesize a path $Y(k), k=1,...,N, N \in \mathbb{N}$, of $Y$, the usual (Euler) method consists in approximating the integral by a Riemann sum. 
Two parameters tune the precision of the method: the discretization 
step $\omega$ and the cut-off value for the integral $\Omega$. The idea in 
\cite{ST3} is to use the fast Fourier transform for an efficient computation
of the Riemann sum. More precisely, let 

\begin{equation*}
 Y(k) = \int_{\mathbb{R}} g(k-s) dM(s) = - \int_{\mathbb{R}} g(s) dM(k-s).
\end{equation*}
Let $\omega,\Omega \in \mathbb{N}$ and 

\begin{equation}\label{repintapprox}
 Y_{\omega,\Omega}(k) = \sum_{j=-\omega\Omega+1}^{0} g(\frac{j-1}{\omega})Z_{\alpha,\omega}(\omega k-j) + \sum_{j=1}^{\omega\Omega} g(\frac{j}{\omega})Z_{\alpha,\omega}(\omega k-j),
\end{equation}
where $Z_{\alpha,\omega}(j) = M(\frac{j+1}{\omega}) - M(\frac{j}{\omega})$ are i.i.d. $\alpha$-stable symmetric random variables. Let $Z_{\alpha}(j)$ denote a sequence of normalised i.i.d $\alpha$-stable symmetric random variables. Then one has the equality in law: $\{ Z_{\alpha,\omega}(j), j \in \mathbb{Z} \}=\{ \omega^{-1/ \alpha}Z_{\alpha}(j), j \in \mathbb{Z}\}$. One may thus write: 

\begin{equation*}
  Y_{\omega,\Omega}(k) = \sum_{j=1}^{2\omega\Omega} a_\omega(j) Z_{\alpha}(\omega(k+\Omega)-j),
\end{equation*}
where 
\begin{displaymath}
 a_{\omega}(j) = 
 \left\{
 \begin{array}{cc}
\omega^{-1/ \alpha} g(\frac{j-1}{\omega}-\Omega) & \textrm{for } j \in \{ 1, ... , \omega\Omega\}\\
\omega^{-1/ \alpha} g(\frac{j}{\omega}-\Omega) & \textrm{for } j \in \{ \omega\Omega+1, ... , 2\omega\Omega\}.\\ 
\end{array}
\right.
\end{displaymath}
For $n \in \bbbz$, let
\begin{displaymath}
 W(n)=\sum_{j=1}^{2\omega\Omega} a_\omega(j) Z_{\alpha}(n-j).
\end{displaymath}
Then $\{Y_{\omega,\Omega}(k), k=1,...,N \}$ has the same law as $\{W(\omega(k+\Omega)), k=1,...,N \}$. 
But $W$ is the convolution product of the sequences $a_\omega$ and $Z_{\alpha}$. As such, it may be be efficiently computed through a fast Fourier transform. See \cite{ST3}
for more details.  

\subsection{Estimation of the approximation error}

When the moving average process is localisable, or more precisely when
the conditions of proposition \ref{prop1.2} are satisfied, it is easy to
assess the performances of the above synthesis method. 

The following proposition gives a bound on the approximation 
error in the $\alpha-$norm. Recall that the $\alpha-$norm 
(defined in (\ref{normdef})) is just
the scale factor of the random variable, and is thus independent of
the integral representation that is used. In addition, it us, up to
a constant depending on $r$ and $\alpha$, an upper bound on moments of
order $0<r<\alpha$.

\begin{prop}\label{errapp}
Let $Y$ be defined by (\ref{ma}), and let $Y_{\omega,\Omega}$ be its approximation defined in (\ref{repintapprox}). Assume $g$ satisfies the conditions of proposition
\ref{prop1.2}. Then, for all $\omega,\Omega \in \mathbb{N}$ and $k \in \mathbb{Z}$ with 
$\omega > \frac{1}{\eta}$, one has 
 \begin{displaymath}
 Err := \| Y(k) - Y_{\omega,\Omega}(k) \|_{\alpha} \leq A_{\omega,\Omega}^{1/ \alpha}
\end{displaymath}
where
\begin{displaymath}
 A_{\omega,\Omega} = \frac{2c^{\alpha}}{(1+a\alpha)\omega^{1+\gamma \alpha}} \sum_{j=1}^{\omega\Omega} \frac{1}{j^{(a-\gamma)\alpha}} + \int_{-\infty}^{-\Omega} |g(s)|^{\alpha} ds + \int_{\Omega}^{+\infty} |g(s)|^{\alpha} ds\\
\end{displaymath}

\end{prop}

\medskip

\noindent {\it Proof.}
By stationarity and independence of the increments of L\'evy motion, one gets:  
\begin{displaymath}
 Err = \sum_{j=-\omega\Omega+1}^{0} \int_{(j-1)/\omega}^{j/\omega} |g(\frac{j-1}{\omega})-g(s)|^{\alpha} ds + \sum_{j=1}^{\omega\Omega} \int_{(j-1)/\omega}^{j/\omega} |g(\frac{j}{\omega})-g(s)|^{\alpha} ds
 \end{displaymath}
\begin{equation}\label{calcdir}
  + \int_{-\infty}^{-\Omega} |g(s)|^{\alpha} ds + \int_{\Omega}^{+\infty} |g(s)|^{\alpha} ds.
\end{equation}

By assumption, for almost all $s \in \mathbb{R}$, $|g(s+h)-g(s)| \leq c|h|^a|s|^{\gamma-a}$ when $0<h<\eta$. Recall that $\omega>\frac{1}{\eta}$. 
A change of variables yields 

\begin{displaymath}
 \int_{(j-1)/\omega}^{j/\omega} |g(\frac{j-1}{\omega})-g(s)|^{\alpha} ds \leq \int_{0}^{1/\omega} c^{\alpha} |s|^{a \alpha} |\frac{j-1}{\omega}|^{(\gamma - a)\alpha} ds
\end{displaymath}
and thus
\begin{eqnarray*}
 Err & \leq & \sum_{j=-\omega\Omega+1}^{0} \int_{0}^{1/\omega} c^{\alpha} |s|^{a \alpha} |\frac{j-1}{\omega}|^{(\gamma - a)\alpha} ds 
 + \sum_{j=1}^{\omega\Omega} \int_{0}^{1/\omega} c^{\alpha} |s|^{a \alpha} |\frac{j}{\omega}|^{(\gamma - a)\alpha} ds\\
 \end{eqnarray*}
\begin{displaymath}
  + \int_{-\infty}^{-\Omega} |g(s)|^{\alpha} ds + \int_{\Omega}^{+\infty} |g(s)|^{\alpha} ds.
\end{displaymath}
Rearranging terms:
\begin{eqnarray*}
 Err & \leq & \frac{c^{\alpha}}{(1+a\alpha)}\frac{1}{\omega^{1+\gamma \alpha}} \left (\sum_{j=-\omega\Omega+1}^{0} |j-1|^{(\gamma-a)\alpha} + \sum_{j=1}^{\omega\Omega} |j|^{(\gamma - a)\alpha}\right) \\ & + & \int_{-\infty}^{-\Omega} |g(s)|^{\alpha} ds + \int_{\Omega}^{+\infty} |g(s)|^{\alpha} ds\\
 & = & A_{\omega,\Omega}.\\
 \end{eqnarray*}
which is the stated result.
\Box

\begin{cor}\label{bound}
Under the conditions of proposition \ref{errapp}, $\| Y(k) - Y_{\omega,\Omega}(k) \|_{\alpha}\to 0$ 
when $(\omega,\Omega)$ tends to infinity.

If in addition $g(x) \leq C |x|^{-\beta}$ when $|x| \to \infty$ for some
$C>0$ and $\beta> \frac {1}{\alpha}$, then:
\begin{equation}\label{errbound}
\| Y(k) - Y_{\omega,\Omega}(k) \|_{\alpha}^{\alpha} \leq K \left(\omega^{-1-\alpha\gamma} +
\Omega^{1-\alpha\beta} \right)
\end{equation}
where $K$ is a constant independent of $k,\omega,\Omega$.
\end{cor}

\noindent {\it Proof.}
Since $g$ satisfies the assumptions of proposition \ref{prop1.2}, $a > \gamma + \frac{1}{\alpha}$. As a consequence, the sum in the first term of $A_{\omega,\Omega}$ converges
when $(\omega,\Omega)$ tends to infinity. The first statement then follows from the facts that
$\alpha\gamma +1>0$ and $g \in \mathcal{F}_{\alpha}$. The second part 
follows by making the obvious estimates. 
\Box

The significance of (\ref{errbound}) is that it allows us to tune
$\omega$ and $\Omega$ to obtain an optimal approximation, 
provided a bound on the
decay of $g$ at infinity is known: optimal pairs $(\omega,\Omega)$ are those
for which the two terms in (\ref{errbound}) are of the same order of magnitude. More precisely,
if the value of $\beta$ is sharp, the order of decay 
of the error will be maximal
when $\Omega=\omega^{\frac{-1-\alpha\gamma}{1-\alpha\beta}}$. Note that the exponent 
$\frac{-1-\alpha\gamma}{1-\alpha\beta}$ is always positive, as expected. Intuitively,
$\omega$ is related to the regularity of $g$ (irregular $g$ requires larger $\omega$), 
while $\Omega$ is linked with the rate of decay of $g$ at infinity. 

For concreteness, let us apply these results to some specific processes:

\begin{exa}\label{Ornrev}(reverse Ornstein-Uhlenbeck process)
Let $Y$ be the reverse Ornstein-Uhlenbeck process defined in proposition
\ref{locOU}. When $\alpha > 1$, we may apply proposition \ref{errapp} 
with $g(x) = \exp(x) \one(x \leq 0), \gamma=0$, $a=1$, $c=2$, $\eta=1$. One
gets, for $\omega>1$, $\Omega>1$,
\begin{displaymath}
 A_{\omega,\Omega} = \frac{2^{\alpha+1}}{1+\alpha} \left( \sum_{j=1}^{\omega\Omega} \frac{1}{j^{\alpha}} \right) \frac{1}{\omega} + \frac{e^{-\alpha \Omega}}{\alpha}.
\end{displaymath}
\end{exa}

However, we may obtain a more precise bound on the approximation
error, valid for any $\alpha \in (0,2)$, 
by using (\ref{calcdir}) directly: 
\begin{displaymath}
 \| Y(k) - Y_{\omega,\Omega}(k) \|_{\alpha}^{\alpha} \leq \frac{2^{\alpha}}{1+\alpha}\left(\frac{1-e^{-\alpha \Omega}}{e^{\alpha /\omega}-1}\right) \frac{1}{\omega^{1+\alpha}}+\frac{e^{-\alpha \Omega}}{\alpha}.
\end{displaymath}
When $(\omega,\Omega) \rightarrow +\infty$, $Err \leq \mathcal{O}(\frac{1}{\omega^{\alpha}})+\mathcal{O}(e^{-\alpha \Omega})$, which is better
than $A_{\omega,\Omega}$ above when $\alpha > 1$. 

We note finally that the optimal choice for $(\omega,\Omega)$ is here $\Omega=\ln(\omega)$, which is
consistent with the fact that the $\beta$ in Corollay \ref{bound} 
may be chosen arbitrarily large.
 
\begin{exa}\label{lfsn}(linear fractional stable noise)
Let $ 0 < \alpha \leq 2$ and let $M$ be an $\alpha$-stable symmetric random measure on $\mathbb{R}$ with control measure $\mathcal{L}$. Let: 
 \begin{displaymath}
  g(x) = (x)_{+}^{H-1/ \alpha} - (x-1)_{+}^{H-1/ \alpha}
 \end{displaymath}
and
\begin{displaymath}
 Y(t) = \int g(t-x) M(dx)  \hspace{0.4cm} (t \in \mathbb{R})
\end{displaymath}
Applying the analysis above with $\gamma=H- \frac{1}{\alpha}$, $a=1$,$c=2|\gamma|$,$\eta=1$ one gets, for $\omega>1$, $\Omega>1$,
\begin{displaymath}
 A_{\omega,\Omega} = \frac{2^{\alpha+1}|H-\frac{1}{\alpha}|^{\alpha}}{1+\alpha}(\sum_{j=1}^{\omega\Omega} \frac{1}{j^{1+\alpha(1-H)}}) \frac{1}{\omega^{\alpha H}} + \int_{\Omega}^{+\infty} |(x)^{H-1/\alpha}-(x-1)^{H-1/\alpha}|^{\alpha} dx
\end{displaymath}
\end{exa}

When $(\omega,\Omega) \rightarrow +\infty$,
\begin{eqnarray*}
 A_{\omega,\Omega} & = & \mathcal{O}(\frac{1}{\omega^{\alpha H}}) + \int_{\Omega}^{+\infty} |(x)^{H-1/\alpha}-(x-1)^{H-1/\alpha}|^{\alpha} dx\\
 & = & \mathcal{O}(\frac{1}{\omega^{\alpha H}}) +  \mathcal{O}(\Omega^{1+\alpha(H-1/\alpha-1)})\\
 & = & \mathcal{O}(\frac{1}{\omega^{\alpha H}}) + \mathcal{O}(\frac{1}{\Omega^{\alpha(1-H)}})\\
\end{eqnarray*}

This process is the one considered in \cite{ST3}. Here we reach a
conclusion similar to \cite[Theorem 2.1]{ST3}, 
which yields the same order of magnitude for the error 
when $(\omega,\Omega) \rightarrow +\infty$. Extensive tests are conducted in \cite{ST3}
to choose the best values for $(\omega,\Omega)$. The criterion for optimizing these
parameters is to test how an estimation method for $H$ performs on synthesized
traces. Here we adopt a different approach based on Corollary \ref{bound}: 
optimal pairs $(\omega,\Omega)$ are those for which (\ref{errbound}) is minimized. Since the value of $\beta = 1-H + 1/\alpha$ is sharp here, one gets 
$\Omega = \omega^{\frac{H}{1-H}}$. It is interesting to
note that the exponent $H/(1-H)$ depends only on the scaling
factor $H$ and not on $\alpha$, and that it may
be larger or smaller than one depending on the value of $H$. We do not have an
explanation for this fact nor for the reason why $H = 1/2$ plays a special r\^ole.

\begin{exa}
As a final illustration, we consider the process of Example \ref{extime}.
With $\gamma=\frac{1}{6}$, $a=1$,$c=1$,$\eta=1$, one gets, for $\omega>1$, $\Omega>1$,
\begin{displaymath}
 A_{\omega,\Omega} = \frac{2}{1+\alpha} \left( \sum_{j=1}^{\omega\Omega} \frac{1}{j^{\frac{5\alpha}{6}}} \right) \frac{1}{\omega^{1+\frac{\alpha}{6}}} + \frac{1}{\frac{5\alpha}{6}-1} \frac{1}{\Omega^{\frac{5\alpha}{6}-1}}.
\end{displaymath}
\end{exa}
\medskip

Again, the value of $\beta=5/6$ is sharp, and the optimal choice is to set
$\Omega=\omega^{\frac{\alpha+6}{5\alpha-6}}$. Since $(\alpha+6)/(5\alpha-6) \geq 1$,
$\Omega$ is larger than $\omega$ in this case, in contrast to the
reverse Ornstein-Uhlenbeck process: it is the decay at infinity of the kernel
that dictates the parameters here, while it was the regularity that mattered in
the case of the reverse Ornstein-Uhlenbeck process.

\subsection{Numerical experiments}

\begin{figure}[!ht]
%\begin{center}
\vspace{-0.15cm}
\includegraphics[height=15.9cm]{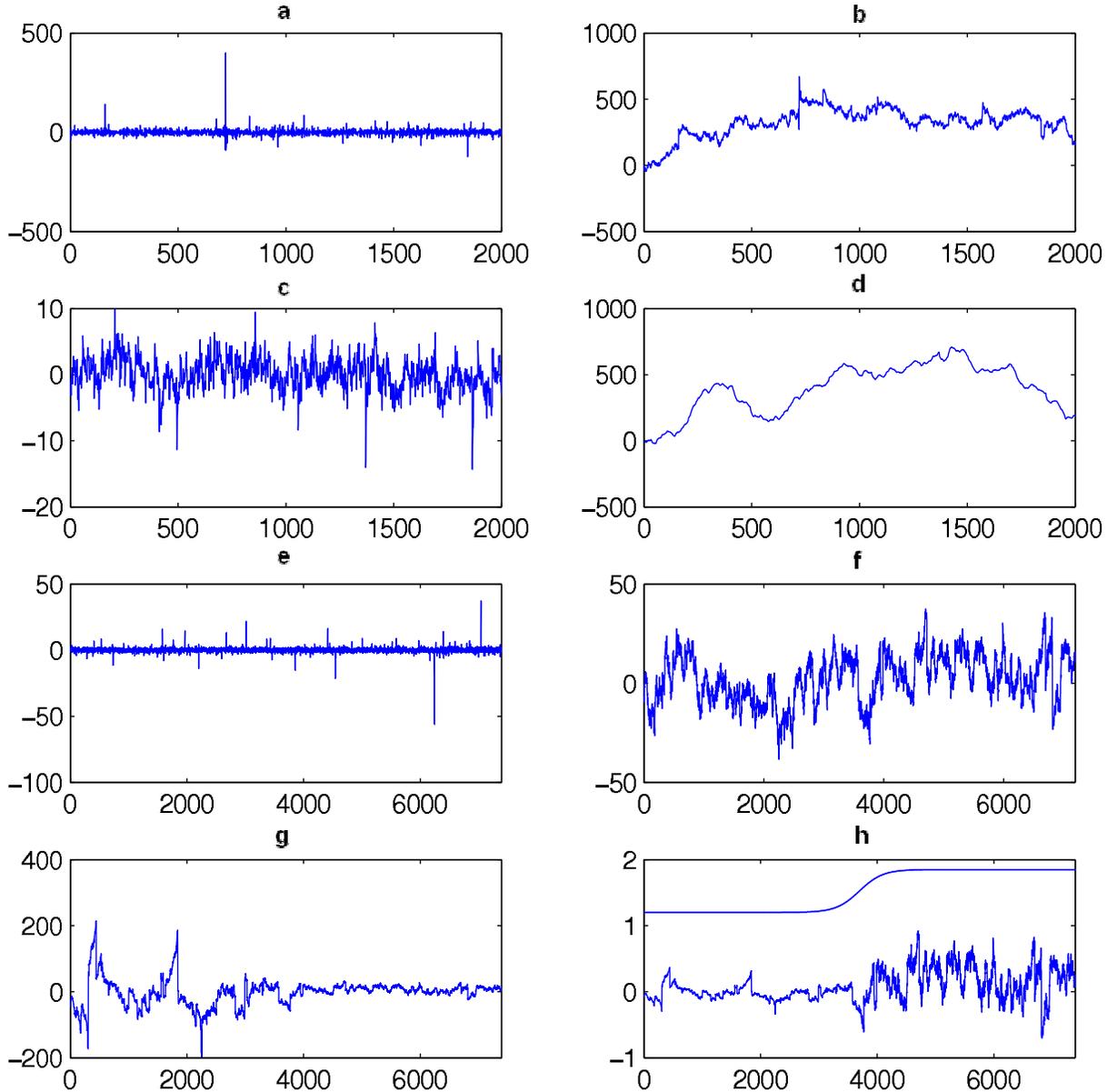}
%\includegraphics[height=22.9cm]{toutalpha2.eps}
%\end{center}
\caption{Paths of localisable processes. (a) The 
process in Example \ref{exfrequency}
and (b) the integrated version. (c) The process 
in Example \ref{extime}
and (d) the integrated version. (e) Reverse Ornstein-Uhlenbeck 
processes with $\lambda=1$, and (f) $\lambda=0.01$.
(g) A multistable reverse Ornstein-Uhlenbeck process with $\lambda=0.01$ 
and (h) the renormalised version along with $\alpha(t)$.}
\label{tout}
\end{figure}

We display in figure \ref{tout} traces of: 
\begin{itemize}
	\item moving average stable processes: the reverse Ornstein-Uhlenbeck 
	process (figures 1(e),(f)), and the processes of Examples \ref{extime} (figure 1(c)) 
	and \ref{exfrequency} (figure 1(a)). In each case, $\alpha = 1.8$.
	Some of the relevant features of the processes of 
Examples \ref{extime} and \ref{exfrequency} seem
to appear more clearly when one integrates them. Thus integral 
versions are displayed in the right-hand part of the corresponding
graphs, figures 1(b),(d).
	\item a multistable version of the reverse Ornstein-Uhlenbeck process,
	using the theory developed
	in Section \ref{ms} (figures 1(g),(h)). Since these processes are localisable, one
	may obtain paths by computing first stable versions with all values
	assumed by $\alpha$, and then ``gluing'' these tangent processes 
	together as appropriate. More precisely, assume we want to obtain the
	values of a multistable process $S$ at the discrete points $(t_1, \ldots, t_n)$. At each $t_i$, $S$ is tangent to a stable process denoted $S_i$. 
	We first synthesize the $n$ stable processes $S_i$ with the method just 
	described, all with the same random seed. The multistable process $S$
	is then obtained by setting $S(t_i) = S_i(t_i), i= 1, \ldots, n$.
Two graphs are displayed for the
multistable process: in figure 1(g)
the graphs are as explained above. In figure 1(h) each ``line'' of the random field
({\it i.e.} the process obtained for a fixed value of $\alpha$)
is renormalized so that it ranges between -1 and 1, prior to
building the multistable process by gluing the paths as
appropriate. This renormalization may be justified using 
Corollary \ref{mapvp}.
\end{itemize}
	
The parameters are as follows:
\begin{itemize}
	\item Process of Example \ref{exfrequency}: $\omega=5000$, $\Omega=877$,
	$N=2000$. The approximation error $Err$ is bounded by 2.172.
	\item Process of Example \ref{extime}: $\omega=104, \Omega=175504$,
	$N=2000$. The term $A_{\omega,\Omega}^{1/\alpha}$ is equal to 0.074.
	\item Reverse Ornstein-Uhlenbeck process with $\lambda=1: \omega=512, \Omega=7,
	N=7392$. The term $A_{\omega,\Omega}^{1/\alpha}$ is equal to 0.0018.
	\item Reverse Ornstein-Uhlenbeck process with $\lambda=0.01: \omega=256, \Omega=800,
	N=7392$. The term $A_{\omega,\Omega}^{1/\alpha}$ is equal to 0.0032. 
	\item Multistable reverse Ornstein-Uhlenbeck process: $\lambda=0.01, \omega=256, \Omega=800,	N=7392$. The $\alpha$ function is the logistic
	function starting from 1.2 and ending at 1.85. More
	precisely, we take:	$\alpha(t) = 1.2 + \frac{0.65}{1+\exp(-\frac{5}{1000}(t-N/2))}$,
	where $N$ is the number of points and $t$ ranges from 1 to $N$ (the graph
	of $\alpha(t)$ is plotted  in figure 1(h).
	Thus, one expects to see large jumps at the beginning of the paths
	and smaller ones at the end. Note that we do not have 
	any results concerning the approximation error for these 
	non-stationary processes. 
\end{itemize}
The value of $\Omega$ in all cases is adjusted so that
the pair $(\omega,\Omega)$ is approximately ``optimal'' as described 
in the preceding  subsection (optimality is not
guaranteed for the multistable processes. 
Nevertheless, since the relation between $\omega$ and $\Omega$ does not
depend on $\alpha$ for the reverse Ornstein-Uhlenbeck process,
it holds in this case). 

The function $g$ of example \ref{exfrequency} cannot be treated using corollary
\ref{bound} nor proposition \ref{errapp} since $g$ does not satisfy 
the conditions of proposition \ref{prop1.2}. However, it is possible to estimate
$Err$ directly. Since $|g(x+h)-g(x)| \leq 2|h|(|x|^{-\frac{3}{2}}\mathbf{1}_{|x|<1}+|x|^{-1}\mathbf{1}_{|x|\geq1})$
and $|g(x)| \leq \frac{4}{|x|}$, one gets:

\begin{displaymath}
 Err^{\alpha} \leq \frac{2^{\alpha+2}}{1+\alpha} \sum_{j=1}^{\omega\Omega} \frac{1}{j^{\alpha}} \frac{1}{\omega^{1-\frac{\alpha}{2}}} + \frac{8}{\alpha-1} \frac{1}{\Omega^{\alpha-1}}
\end{displaymath}
\begin{displaymath}
  Err^{\alpha} \leq \frac{2^{\alpha+2}}{1+\alpha} \frac{\alpha}{\alpha-1}  \frac{1}{\omega^{1-\frac{\alpha}{2}}} + \frac{8}{\alpha-1} \frac{1}{\Omega^{\alpha-1}}.
\end{displaymath}
The asymptotic optimal relation between $\omega$ and $\Omega$ is thus $\Omega=\omega^{\frac{2-\alpha}{2(\alpha-1)}}=\omega^{0.125}$. The values
in our simulation are slightly different since they are chosen
to optimize the actual expression with a finite $\omega$.

Finally, we stress that the same random
seed ({\it i.e.} the same underlying stable $M(dx)$) has been used 
for all simulations, for easy comparison. Thus, for instance, the jumps 
appear at precisely the same locations in each graph. Notice in
particular the ranges assumed by the different processes.

The differences between the graphs of the processes of
Examples \ref{extime}, \ref{exfrequency} and the 
reverse Ornstein-Uhlenbeck process are easily 
interpreted by examining the 
three kernels: the kernel of the process of Example \ref{exfrequency} 
diverges at 0, thus putting more emphasis on strong jumps, as seen 
on the picture, with more jaggy curves and an ``antipersistent'' 
behaviour. The kernel of the process of Example \ref{extime}, in contrast,
is smooth at the origin. In addition, it has a slow decay. These features 
result in an overall smoother appearance and allow ``trends'' to appear 
in the paths. Finally, the kernel of the reverse Ornstein-Uhlenbeck
process has a decay controlled by $\lambda$. For ``large''
$\lambda$ (here, $\lambda=1$), little averaging is done, and the resulting
path is very irregular. For ``small'' $\lambda$ (here, $\lambda=0.01$),
the kernel decays slowly and the paths look smoother (recall that, in
the Gaussian case, the Ornstein-Uhlenbeck tends in distribution to white
noise when $\lambda$ tends to infinity, and to Brownian motion when
$\lambda$ tends to 0).

\end{document}